\documentclass{ifacconf}

\usepackage[dvipsnames]{xcolor}
\usepackage{listings}
\usepackage{graphicx}      
\usepackage{natbib}        
\usepackage{amsfonts,amsmath,amssymb}
\usepackage{algorithmic}
\usepackage{url,cite}
\usepackage{graphicx}
\usepackage{tikz,pgfplots}
\usepackage{hhline}

\newcommand{\bx}{\boldsymbol{x}}
\newcommand{\bM}{\boldsymbol{M}}

\newcommand{\blambda}{\boldsymbol{\lambda}}
\newcommand{\bz}{\boldsymbol{z}}
\newcommand{\bX}{\boldsymbol{X}}

\newcommand{\bZ}{\boldsymbol{Z}}

\newcommand{\bd}{\boldsymbol{d}}
\newcommand{\bp}{\boldsymbol{p}}
\newcommand{\bff}{\boldsymbol{f}}
\newcommand{\bc}{\boldsymbol{c}}
\newcommand{\bg}{\boldsymbol{g}}

\newcommand{\bn}{\boldsymbol{n}}
\newcommand{\bm}{\boldsymbol{m}}
\newcommand{\br}{\boldsymbol{r}}

\newcommand{\st}{\mathop{\text{\normalfont s.t.}}}
\newcommand{\diag}{\mathop{\text{\normalfont diag}}}

\graphicspath{{fig}}
\allowbreak
\allowdisplaybreaks

\lstdefinelanguage{Julia}%
{morekeywords={abstract,break,case,catch,const,continue,do,else,elseif,%
end,export,false,for,function,immutable,import,importall,if,in,%
macro,module,otherwise,quote,return,switch,true,try,type,typealias,%
using,while,@threads,@optinode,@optiedge},%
sensitive=true,%
alsoother={$},
morecomment=[l]\#,%
morecomment=[n]{\#=}{=\#},%
morestring=[s]{"}{"},%
morestring=[m]{'}{'},%
}[keywords,comments,strings]%
\definecolor{codegreen}{rgb}{0,0.6,0}
\definecolor{codegray}{rgb}{0.5,0.5,0.5}
\definecolor{codepurple}{rgb}{0.58,0,0.82}
\definecolor{backcolour}{rgb}{0.9,0.9,0.9}
\lstdefinestyle{mystyle}{
backgroundcolor=\color{backcolour},   
commentstyle=\color{codegreen},
keywordstyle=\color{magenta},
numberstyle=\tiny\color{codegray},
stringstyle=\color{codepurple},
basicstyle=\ttfamily\footnotesize,
breakatwhitespace=false,         
breaklines=true,                 
captionpos=b,                    
keepspaces=true,                 
numbers=left,                    
numbersep=5pt,                  
showspaces=false,                
showstringspaces=false,
showtabs=false,                  
tabsize=2
}
\begin{document}
\begin{frontmatter}

\title{Graph-Based Modeling and Decomposition of Energy Infrastructures\thanksref{footnoteinfo}} 

\thanks[footnoteinfo]{We acknowledge support from the Grainger Wisconsin Distinguished Graduate Fellowship.}

\author[Madison]{Sungho Shin} 
\author[LANL]{Carleton Coffrin} 
\author[LANL]{Kaarthik Sundar}
\author[Madison]{Victor M. Zavala}

\address[Madison]{University of Wisconsin-Madison, Madison, WI 53706 USA\\ (e-mail: sungho.shin@wisc.edu; victor.zavala@wisc.edu).}
\address[LANL]{Los Alamos National Laboratory, Los Alamos, NM 87545 USA\\ (e-mail: cjc@lanl.gov; kaarthik@lanl.gov)}

\begin{abstract}                
Nonlinear optimization problems are found at the heart of real-time operations of critical infrastructures.  These problems are computationally challenging because they embed complex physical models that exhibit space-time dynamics. We propose modeling these problems as graph-structured optimization problems, and illustrate how their structure can be exploited at the modeling level (for parallelizing function/derivative computations) and at the solver level (for parallelizing linear algebra operations). Specifically, we present a restricted additive Schwarz scheme that enables flexible decomposition of complex graph structures within an interior-point algorithm. The proposed approach is implemented as a general-purpose nonlinear programming solver that we call MadNLP.jl; this Julia-based solver is interfaced to the graph-based modeling package Plasmo.jl.  The efficiency of this framework is demonstrated via problems arising in transient gas network optimization and multi-period AC optimal power flow. We show that our framework accelerates the solution (compared to off-the-shelf tools) by over 300\%; specifically, solution times are reduced from 72.36 sec to 23.84 sec for the gas problem and from 515.81 sec to 149.45 sec for the power flow problem. 
\end{abstract}


\end{frontmatter}

\section{Introduction}
Real-time operation of modern energy infrastructures requires solving large-scale nonlinear programs (NLPs). Application examples include transient gas network optimization  \citep{sundar2018state} and multi-period optimal power flow problems \citep{geth2020flexible,kim2020real}. Achieving real-time solutions for these problems is challenging, as they embed complex physical models that require space-time discretization. NLPs arising in this context can easily reach millions of variables and constraints and defy the scope of off-the-shelf solvers. Specifically, scalability bottlenecks are often encountered at the modeling level (function and derivative computations) and at the solver level (computation of the search step).

Large-scale NLPs arising in energy infrastructures have the key characteristic that they exhibit {\em sparse graphs} structures; we refer to such problems as {\it graph-structured optimization problems} \citep{jalving2019graph,shin2020decentralized,jalving2020graph}. Graph-structured problems can be conveniently modeled using specialized modeling platforms such as {\tt Plasmo.jl} \citep{jalving2019graph,jalving2020graph} and solved using structure-exploiting optimization solvers such as {\tt PIPS-NLP} \citep{chiang2014structured}). {\tt Plasmo.jl} is a graph-based modeling platform that enables the  modular construction and analysis of highly complex models; this platform also leverages the algebraic modeling capabilities of {\tt JuMP.jl} \citep{dunning2017jump} and  facilitates access to infrastructure modeling tools such as {\tt GasModels.jl} and {\tt PowerModels.jl} \citep{russell2020gasmodels,8442948}.

Another key benefit of {\tt Plasmo.jl} is that it can communicate model structures to solvers, and this facilitates the implementation of different decomposition strategies, notably the alternating direction method of multipliers \citep{boyd2011distributed}, overlapping Schwarz \citep{shin2020overlapping}, and parallel interior-point (IP) methods \citep{chiang2014structured,rodriguez2020scalable}. 

In this paper, we present a new and flexible decomposition framework for graph-structured optimization problems. Our framework uses a restricted additive Schwarz (RAS) decomposition scheme implemented within a filter line-search IP method \citep{wachter2006implementation}. We present a {\em Julia}-based implementation of this approach, which we call {\tt MadNLP.jl} (\url{https://github.com/zavalab/MadNLP.jl}).  We use our framework to experiment with different decomposition strategies that exploit parallelism at the modeling and at the solver level. Specifically, we consider a scheme that parallelizes function and derivative computations by exploiting the modular structure of {\tt Plasmo.jl}. We also consider a scheme that uses RAS \citep{cai1999restricted} for parallelizing step computations. RAS has been widely used for the solution of large linear algebra systems arising from discretized partial differential equations (PDEs) \citep{balay2019petsc}, but we have recently shown that it can also be applied to solve general linear systems arising in graph-structured optimization  \citep{shin2020decentralized,gerstner2016domain}.  Our computational results indicate that our framework can accelerate computations by up to 300\% (compared to off-the-shelf tools). 

The paper is organized as follows: In Section \ref{sec:modeling}, we define the graph-structured problem of interest and discuss how its structure can be exploited at the modeling level. In Section \ref{sec:decomposition}, we discuss parallel decomposition schemes. In Section \ref{sec:results}, we apply these schemes to transient gas network optimization and multi-period AC optimal power flow problems. Section \ref{sec:conclusions} presents concluding remarks.

\begin{figure}[!htp]
  \centering
  \begin{tikzpicture}
    \node[xscale=-1] at (0,0) {\includegraphics[width=1in]{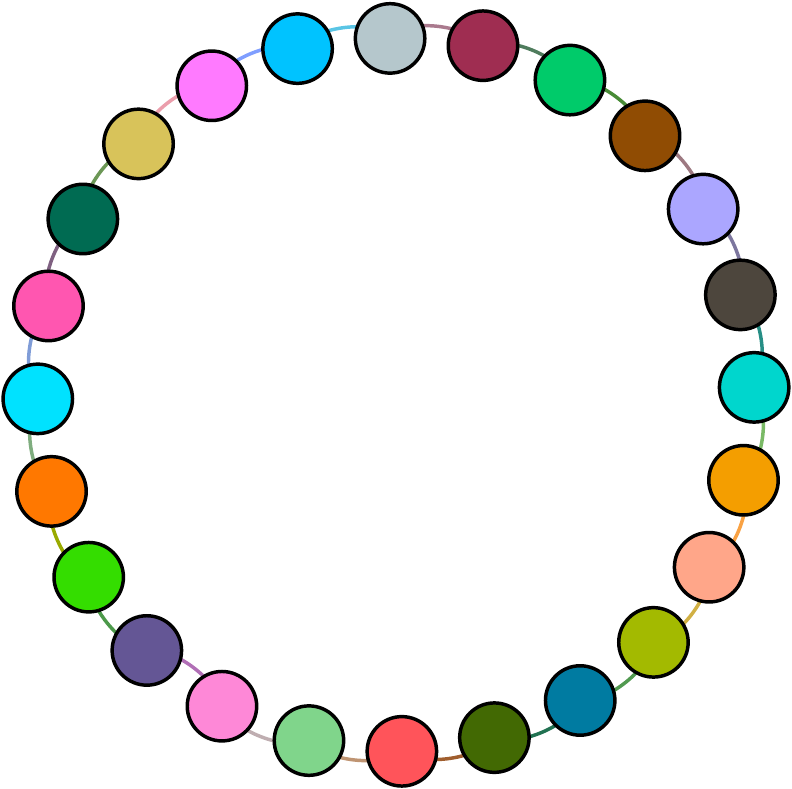}};
    \node at (6,0) {\includegraphics[width=1in]{fig/adchem_lower.pdf}};
    \node at (3,0) {\includegraphics[width=1.2in]{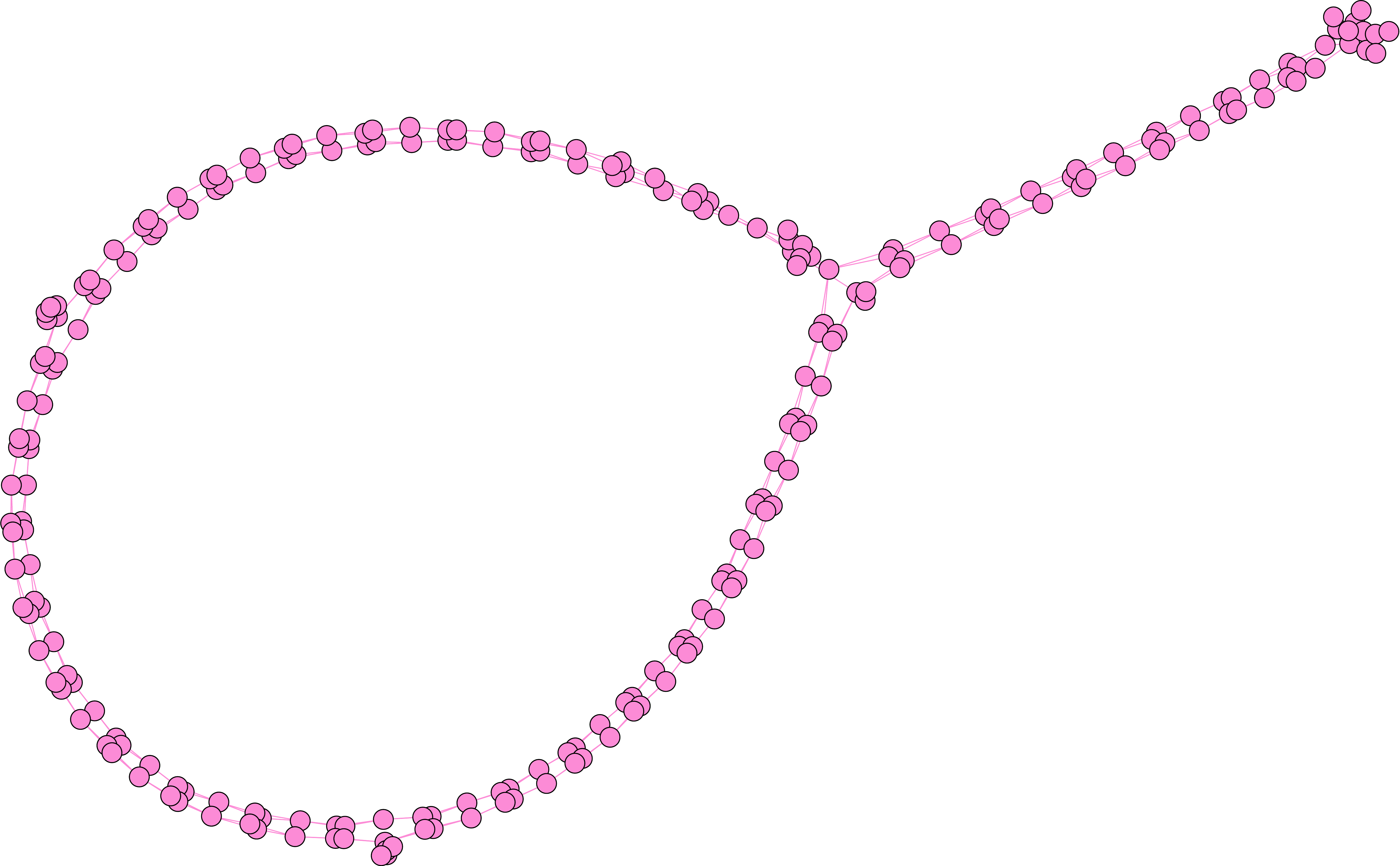}};
    \draw[<-,dashed] (1.8,.6) to [out=150,in=30] (.7,1.1);
    \draw[->,dashed] (5.7,.9) to [out=120,in=30] (3.5,.6);
  \end{tikzpicture}
  \caption{Alternative graph representations of a transient gas network problem. Left: graph $G$  in which nodes represent times and the spatial (network) structure is embedded within each node. Right: graph $\widetilde{G}$ in which each node represents a variable or a constraint. Middle: a slice of $\widetilde{G}$ that corresponds to the structure contained in each node of $G$.}\label{fig:g-gtilde}
\end{figure}

\section{Graph-Based Modeling}\label{sec:modeling}

Optimization problems arising in energy infrastructures can be expressed as a graph-structured optimization problem of the form:
\begin{subequations}\label{eqn:prob}
  \begin{align}
    \min_{\{x_i\}_{i\in {V}}}\;
    & \sum_{i\in {V}}f_i(x_i)\label{eqn:prob-obj}\\
    \st\;
    & c^I_i(x_i) = 0,\;i\in {V},\;(\lambda^I_i)\label{eqn:prob-inner}\\
    & c^L_i(\{x_j\}_{j\in N_{G}[i]}) = 0,\;i\in {V},\;(\lambda^L_i)\label{eqn:prob-eq}\\
    &x_i\geq 0,\; i\in {V}\;(z_i)\label{eqn:prob-bound}
\end{align}
\end{subequations}
Here, the undirected graph ${G}=({V},{E})$ is an ordered pair of the nonempty, strictly ordered node set ${V}$ and the edge set ${E}\subseteq \{\{i,j\}\subseteq {V}: i\neq j\}$; $N_{G}[i]:=\{j\in {V}:\{i,j\}\in {E}\}\cup\{i\}$ denotes the neighborhood of $i\in {V}$ on ${G}$. For each $i\in {V}$, $x_i\in\mathbb{R}^{n_i}$ is the decision variable; $f_i(\cdot)$ is the objective function; $c^I_i(\cdot)$ is the inner equality constraint function; $c^L_i(\cdot)$ is the linking constraint function; $\lambda^I_i\in\mathbb{R}^{m^I_i}$ is the dual variable associated with \eqref{eqn:prob-inner}; $\lambda^L_i\in\mathbb{R}^{m^L_i}$ is the dual variable associated with \eqref{eqn:prob-eq}; $z_i\in\mathbb{R}^{n_i}$ is the dual variable associated with \eqref{eqn:prob-bound}.  General inequality constraints can be handled by introducing slack variables. 

In the context of energy infrastructures, the graph ${G}$ may intuitively be used as an abstraction of the space-time structure of the problem. Specifically, each node $i\in V$ is a component located at a particular spatial location and at a particular time point (e.g., a subset of buses, generators, storage facilities, and electric lines in power networks or a subset of junctions, compressors, and pipelines in gas networks). The link constraints \eqref{eqn:prob-eq} may represent spatial connections (e.g., interconnecting electric lines in power networks or pipelines in gas networks) or temporal connections. Furthermore, we will show that the graph ${G}$ can also be used as a general abstraction in which each node $i\in V$ represents an individual variable or constraint of the problem (the graph encodes the general sparsity structure of the problem). The ability to represent the same optimization problem in different forms provides flexibility to identify efficient decomposition strategies. 

We define the short-hand notation for \eqref{eqn:prob}:
\begin{subequations}\label{eqn:prob-simple}
\begin{align}
\min_{\bx}\; &\bff(\bx)\\  \st\; &\bc(\bx) = 0\quad (\blambda)\\ &\bx\geq 0\quad(\bz),
\end{align}
\end{subequations}
where $\bx:=\{x_i\}_{i\in V}$, $\blambda=\{[\lambda^I_i;\lambda^L_i]\}_{i\in V}$, $\bff(\cdot):=\sum_{i\in V}f_i(\cdot)$, $\bc(\cdot):=\{[c^I_i(\cdot);c^L_i(\cdot)]\}_{i\in V}$, $\bn=\sum_{i\in V} n_i$, and $\bm=\sum_{i\in V} m^I_i+m^L_i$. Here, $n_i$ and $m_i$ are the primal-dual variable dimensions, and $\{(\cdot)_i\}_{i\in V}$ denotes the vector concatenation. We use boldface symbols to denote quantities associated with multiple nodes. 

Each node $i\in V$ may contain more than or less than one variable and constraint ($n_i,m_i\in\mathbb{Z}_{\geq 0}$). Thus, the {\it problem graph} $G=(V,E)$ may be different from the {\it primal-dual coupling graph} $\widetilde{G}=(\widetilde{V},\widetilde{E})$, where $\widetilde{V}:=\mathbb{Z}_{[1,\bn+\bm]}$, $\widetilde{E}:=\{\{i,j\}: \nabla^2_{(\bx,\blambda)(\bx,\blambda)}\mathcal{L}(\bx,\blambda,\bz)[i,j]\neq 0\}$. Here, $\mathcal{L}(\cdot)$ is the Lagrangian of \eqref{eqn:prob-simple}, and we use syntax $\mathbb{Z}_{[a,b]}:=\{a,a+1,\cdots,b\}$. We can observe that a node $i$ in $V$ corresponds to a set of nodes $U_i\subseteq \mathbb{Z}_{[1,\bn+\bm]}$ in $\widetilde{V}$, which contains multiple variables and constraints. Example graphs $G$ and $\widetilde{G}$ for a transient gas network problem are depicted in Figure \ref{fig:g-gtilde} (a detailed problem formulation can be found in Section \ref{sec:gas}). Graph $G$ contains 24 nodes, each corresponding to a time point in a prediction horizon. Periodicity (over a 24 hours period) is enforced as constraints (this periodicity creates the cycle shape of $G$). In this graph, each node embeds the spatial structure of the problem (network and pipelines). Graph $\widetilde{G}$ unfolds the temporal and spatial structure and shows the interconnectivity between all variables and constraints in the problem.

\section{Graph-Based Decomposition}\label{sec:decomposition}
We proceed to describe our IP solver {\tt MadNLP.jl}, comprising a new {\tt Schwarz} decomposition scheme that exploits graph structures within an IP method, and we describe its interface to {\tt Plasmo.jl}.

\subsection{Interior-Point Method}
The IP method implemented in {\tt MadNLP.jl} finds the solution of \eqref{eqn:prob-simple} by solving a sequence of {\it barrier subproblems}:
\begin{align}\label{eqn:barrier}
  \min\;&\varphi(\bx):=\bff(\bx) - \mu\,\mathbf{e}^T\log(\bx)\quad \st\; \bc(\bx) = 0.
\end{align}
with a decreasing sequence for parameter $\mu$. The KKT conditions for \eqref{eqn:barrier} give the nonlinear equations:
\begin{align}\label{eqn:kkt}
  \nabla\bff (\bx) + A^\top \blambda -\bz =0;\; \bc(\bx)  = 0;\; \bX\bZ \mathbf{e} - \mu\, \mathbf{e}= 0,
\end{align}
where $A:=\nabla\bc(\bx)$, $\bX:=\diag(\bx)$, and $\bZ:=\diag(\bz)$. A solution of KKT system  \eqref{eqn:kkt} is obtained by computing primal-dual Newton steps $\bd^*$ from:
\begin{align}\label{eqn:newton}
  \underbrace{\begin{bmatrix}
    W+\Sigma+\delta_w I &A^\top\\
    A& -\delta_c I
  \end{bmatrix}}_{\bM}
  \underbrace{\begin{bmatrix}
    d^x\\d^\lambda
  \end{bmatrix}}_{\bd}
  =-
  \underbrace{\begin{bmatrix}
      \nabla_{\bx} \varphi(\bx) + A^\top \blambda \\
    \bc(\bx).
  \end{bmatrix}}_{\bp},
\end{align}
where $W:=\nabla^2_{\bx\bx}\mathcal{L}(\bx,\blambda,\bz)$, $\Sigma:=\bX^{-1}\bZ$, and $\delta_w,\delta_c>0$ are regularization parameters. The step $\bd^*$ computed from \eqref{eqn:newton} is safeguarded by a line-search filter procedure to induce global convergence  \citep{wachter2006implementation}. Typically, the solution of the linear system \eqref{eqn:newton} is the most computationally intensive step in the IP method. This system is typically solved using direct linear solvers that are based on $LDL^\top$ factorizations (e.g., as those implemented in HSL routines \citep{hsl2007collection}). Decomposition strategies based on Schur complements   \citep{chiang2014structured} and iterative strategies \citep{curtis2012note,rodriguez2020scalable} have also been proposed. The solution of \eqref{eqn:newton} based on a direct block $LDL^\top$ factorization reveals the inertia (the number of positive, zero, negative eigenvalues) of $\bM$. This inertia information is crucial in determining the acceptability of the computed step and in triggering the regularization of the linear system. However, inertia is not available when using iterative solution algorithms (as that proposed in this work). In {\tt MadNLP.jl}, we use an inertia-free regularization strategy to determine the acceptability of the step \citep{chiang2016inertia}. This method performs a simple negative curvature test to trigger regularization. 

\subsection{Restricted Additive Schwarz (RAS)}
We propose a solution strategy for \eqref{eqn:newton} based on an RAS scheme. We define some concepts and quantities that help explain our algorithm. Consider a partition $\{V_k\}_{k=1}^K$ of $V$; we call $V_k$ {\it non-overlapping subdomains}. This partition can be obtained by applying a graph partitioning scheme to $G$. We then construct a family of {\it overlapping subdomains}  $\{V_k^\omega\}_{k=1}^K$ (these are constructed by {\it expanding} $V_k$). The expansion procedure is performed by progressively incorporating adjacent nodes \citep{shin2020decentralized} (the {\it size of overlap} $\omega$ represents the expansion level). We observe that: 
\begin{align*}
  V_k\subseteq  V^\omega_k \subseteq V,\; k=1,\cdots, K;\quad \bigsqcup_{k=1}^K V_k=\bigcup_{k=1}^K V^\omega_k = V,
\end{align*}
where $\sqcup$ denotes disjoint union. With $\{V_k\}_{k=1}^K$ and $\{V^\omega_k\}_{k=1}^K$, we define the corresponding index sets in the space of primal-dual variables in $\mathbb{R}^{\bn+\bm}$ as follows: 
\begin{align*}
  W_k := \bigsqcup_{i\in V_k} U_i;\quad W^\omega_k := \bigsqcup_{i\in V^\omega_k} U_i,\;k=1,\cdots,K,
\end{align*}
where $U_i\subseteq W:=\mathbb{Z}_{[1,\bn+\bm]}$ is the index set of $[x_i;\lambda_i]$ in $[\bx;\blambda]$. We now observe that:
\begin{align*}
  W_k\subseteq  W^\omega_k \subseteq W,\; k=1,\cdots, K;\;\quad \bigsqcup_{k=1}^K W_k=\bigcup_{k=1}^K W^\omega_k = W.
\end{align*}

We state the RAS scheme for solving \eqref{eqn:newton} as:
\begin{align}\label{eqn:ras}
  \bd^{(\ell+1)} = \bd^{(\ell)} + \Big(\underbrace{\sum_{k=1}^K \widetilde{R}_k \bM^{-1}_k R_k^\top}_{P^{-1}}\Big) \br^{(\ell)},\; \ell=0,1,\cdots.
\end{align}
Here, $\ell$ is the RAS iteration counter, $\br^{(\ell)} := \bp-\bM\bd^{(\ell)}$ is the residual; $\bM_k := \bM [W^\omega_k,W^\omega_k]$; $R_k := \{e^\top_i\}^\top_{i\in W_k}$; $\widetilde{R}_k:= \{(\widetilde{e}^k_i)^\top\}^\top_{i\in W_k}$, where $e_i$ is the $i$-th standard basis of $\mathbb{R}^{\bn+\bm}$, and $\widetilde{e}^k_i=e_i$ if $i\in W_k$ and $0$ otherwise. 

The RAS scheme \eqref{eqn:ras} involves the following steps. We first obtain the residual at the current step $\ell$; then, for each {\it overlapping} subdomains $\{V^\omega_k\}_{k=1}^K$, the associated residual is extracted as $R^\top_k \br^{(\ell)}$. The $k$-th subsystem is then solved by applying $\bM^{-1}_k$. In {\tt MadNLP.jl}, a factorization of $\bM_k$ is computed with a direct solver and stored, so that the system can be repeatedly solved whenever the new right-hand-side is given. Subsequently, the solution $\bM^{-1}_k R^\top_k \br^{(\ell)}$ for the $k$-th overlapping subdomain is {\it restricted} to the non-overlapping subdomain $V_k$ and then mapped back to the full-space by applying $\widetilde{R}^\omega_k$ (the indices associated with $V\setminus V_k$ are set to zero in this step).
Key defining features of RAS are the concepts of {\it overlap} and {\it restriction}. The overlap allows the dampening of the adverse effect of the truncated domain, and the restriction procedure discards the part of the solution where the adverse truncation effect is strong.

It has been empirically observed that the convergence of the RAS algorithm improves as the size of overlap $\omega$ increases and as the conditioning of $\bM$ improves \citep{cai1996overlapping}. For positive definite $\bM$, an exponential relationship between the convergence rate and the size of overlap has been established \citep{shin2020decentralized}. When $\omega=0$, the RAS scheme reduces to a block-Jacobi scheme (decentralized) while, when $\omega$ is maximal ($\bM_k=\bM$), the RAS scheme becomes a direct solution method (centralized). In this sense, RAS provides a bridge between fully decentralized and fully centralized schemes (thus providing flexibility). In the {\tt Schwarz} submodule of ${\tt MadNLP.jl}$, $\omega$ is set automatically based on the relative size of $V_k$, and adaptively adjusted whenever a convergence issue occurs. Algorithm \eqref{eqn:ras} uses a simple static iteration (also called a Richardson iteration), but more sophisticated iterative methods (e.g., the generalized mean residual (GMRES) method) can also be used (by treating $P$ as a preconditioner). Both Richardson and GMRES iterators are implemented in {\tt MadNLP.jl}.

\begin{figure}[t]
\centering\footnotesize
\begin{tikzpicture}[scale=.5]
  \node[xscale=-1] at (1.2,0) {\includegraphics[width=.7in]{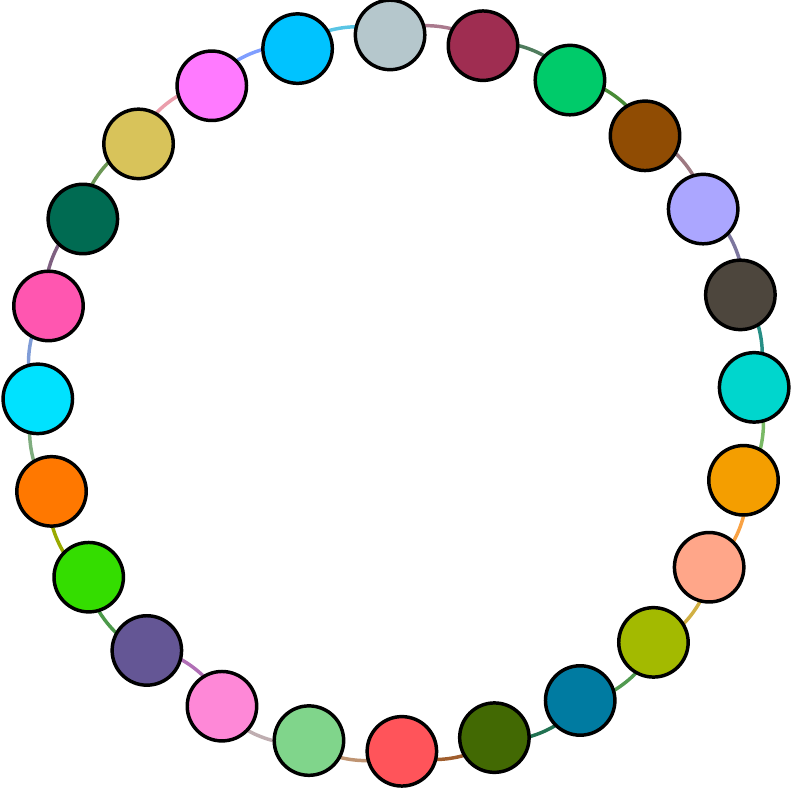}};
  \node at (1.2,0) {$G=(V,E)$};
  
  \node[xscale=-1] at (5.6,0) {\includegraphics[width=.7in]{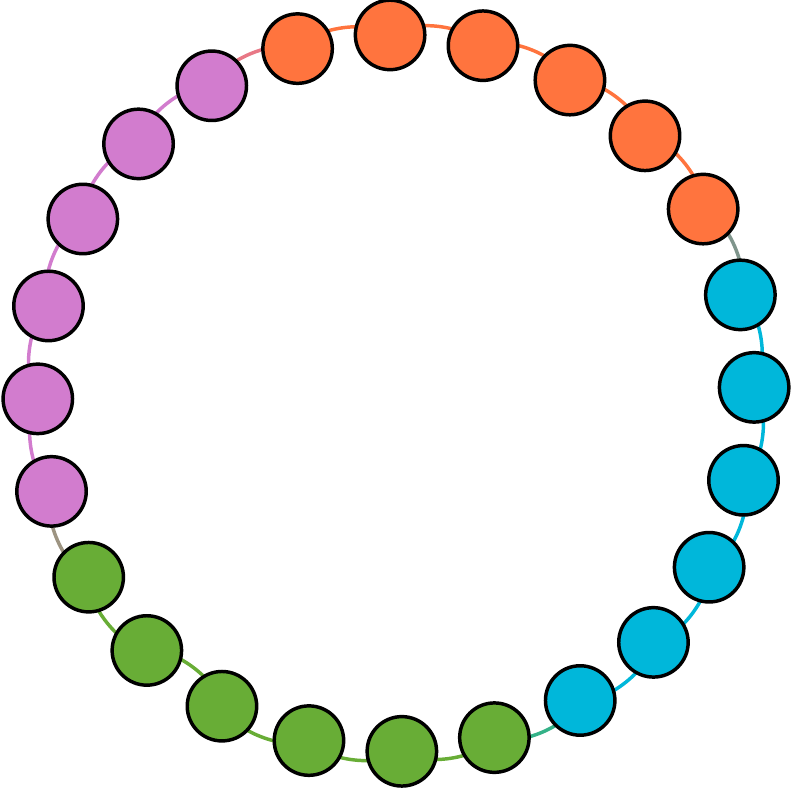}};
  \node at (5.2,1) {\color{RedOrange}$V_1$};
  \node at (6.5,.5) {\color{Purple}$V_2$};
  \node at (6,-1) {\color{LimeGreen}$V_3$};
  \node at (4.7,-.5) {\color{Cyan}$V_4$};
  \node at (10,0) {\includegraphics[width=.8in]{fig/adchem_lower_part.pdf}};
  \node at (8,1.5) {\color{RedOrange}$W_1$};
  \node at (11.5,1.5) {\color{Purple}$W_2$};
  \node at (12.1,-1.3) {\color{LimeGreen}$W_3$};
  \node at (8,-1.3) {\color{Cyan}$W_4$};
  \draw (12.6,1.8) rectangle (16.2,-1.8);
  \foreach \x in {0,...,22}
  \fill (\x*.15+12.6,1.8-\x*.15) rectangle (\x*.15+.3+12.6,1.8-\x*.15-.3);
  \fill[RedOrange,opacity=.5] (12.6,1.8) rectangle (13.8,.6);
  \fill[Purple,opacity=.5] (13.2,1.2) rectangle (14.7,-.3);
  \fill[LimeGreen,opacity=.5] (14.1,.3) rectangle (15.6,-1.2);
  \fill[Cyan,opacity=.5] (15,-.6) rectangle (16.2,-1.8);
  \fill (12.6,-1.8) rectangle (12.75,-1.65);
  \fill (16.05,1.65) rectangle (16.2,1.8);
  \fill[RedOrange,opacity=.5] (12.6,-1.8) rectangle (13.8,-1.5);
  \fill[RedOrange,opacity=.5] (15.9,.6) rectangle (16.2,1.8);
  \fill[Cyan,opacity=.5] (12.6,-1.8) rectangle (12.9,-.6);
  \fill[Cyan,opacity=.5] (15,1.5) rectangle (16.2,1.8);
  \node at (14.2,1.5) {\color{RedOrange}$\bM_1$};
  \node at (15.3,.7) {\color{Purple}$\bM_2$};
  \node at (13.5,-.7) {\color{LimeGreen}$\bM_3$};
  \node at (14.4,-1.5) {\color{Cyan}$\bM_4$};
\end{tikzpicture}
\caption{Graph-based decomposition in {\tt MadNLP.jl} applied to a transient gas network problem.}\label{fig:decomposition}
\end{figure}

The construction of partitions for the RAS scheme is illustrated in Figure \ref{fig:decomposition}; these partitions correspond to the transient gas network example of Figure \ref{fig:g-gtilde}. Here, by partitioning the problem graph $G=(V,E)$ (first subfigure), the node set $V$ is divided into 4 subdomains $V_1,V_2,V_3,V_4$ (second subfigure). These subdomains are associated with the primal-dual index sets $W_1,W_2,W_3,W_4$ (third subfigure). After applying expansions, the associated blocks $\bM_1,\bM_2,\bM_3,\bM_4$ of $\bM$ are identified (last subfigure); these blocks are used by the RAS scheme \eqref{eqn:ras}.

\subsection{Implementation in {\tt MadNLP.jl}}
\begin{figure*}[t]
  \centering\small
\begin{tikzpicture}
    \def\x{3.6};
    \def\lx{1.7};
    \def\ly{.65};
    \def\s{0.4};
    \def\w{0.05};
    \fill[color=backcolour] (-\lx-\x,-\ly) rectangle (\lx-\x,\ly);
    \fill[color=backcolour] (-\lx,-\ly) rectangle (\lx,\ly);
    \fill[color=backcolour] (-\lx+\x,-\ly) rectangle (\lx+\x,\ly);
    \node at (-\x,0) {\tt Plasmo.OptiGraph};
    \node at (0,0) {\tt MadNLP.Solver};
    \node at (\x,0) {\tt Schwarz.Solver};
    \foreach \i in {0,2}{
      \fill[color=backcolour] (-\lx-2*\x,-\ly+\i*\s+\i*\w) rectangle (\lx-2*\x,-\ly+\i*\s+\s+\i*\w);
      \fill[color=backcolour] (-\lx+2*\x,-\ly+\i*\s+\i*\w) rectangle (\lx+2*\x,-\ly+\i*\s+\s+\i*\w);
      \node at (-2*\x,-\lx+\i*\s+\s/2+\i*\w+1.05) {\tt JuMP.Model};
      \node at (2*\x,-\lx+\i*\s+\s/2+\i*\w+1.05) {\tt Ma57.Solver};
    }
    \node at (-2*\x,0.1) {$\vdots$};
    \node at (2*\x,0.1) {$\vdots$};
    
    \draw[<-] (-2*\x+\s,-\ly-2*\w)-- node[below] {$\{\bx^*_k,\blambda^*_k,\bz^*_k\}$} (-\x-\s,-\ly-2*\w);
    \draw[<-] (-\x+\s,-\ly-2*\w)-- node[below] {$\bx^*,\blambda^*,\bz^*$} (-\s,-\ly-2*\w);
    \draw[<-] (+\s,-\ly-2*\w)-- node[below] {$\bd^*$} (\x-\s,-\ly-2*\w);
    \draw[<-] (\x+\s,-\ly-2*\w)--node[below] {$\{\bd^*_k\}$} (2*\x-\s,-\ly-2*\w);

    \draw[->] (-2*\x+\s,\ly+2*\w)-- node[above] {$\{\bff_k(\cdot),\bg_k(\cdot)\}$} (-\x-\s,\ly+2*\w);
    \draw[->] (-\x+\s,\ly+2*\w)-- node[above] {$\bff(\cdot),\bg(\cdot)$} (-\s,\ly+2*\w);
    \draw[->] (+\s,\ly+2*\w)-- node[above] {$\bM,\bp$} (\x-\s,\ly+2*\w);
    \draw[->] (\x+\s,\ly+2*\w)-- node[above] {$\{\bM_k,\br^{(\ell)}_k\}$} (2*\x-\s,\ly+2*\w);
  \end{tikzpicture}
  \centering
  \vspace{0.1in}
  \begin{tikzpicture}
    \def\x{3.6};
    \def\lx{1.7};
    \def\ly{.45};
    \def\s{0.75};
    \def\w{0.05};
    \fill[color=backcolour] (-\lx-2*\x,-\ly) rectangle (\lx-2*\x,\ly);
    \fill[color=backcolour] (-\lx,-\ly) rectangle (\lx,\ly);
    \fill[color=backcolour] (-\lx+2*\x,-\ly) rectangle (\lx+2*\x,\ly);
    \node at (-2*\x,0) {\tt JuMP.Model};
    \node at (0,0) {\tt MadNLP.Solver};
    \node at (2*\x,0) {\tt Ma57.Solver};    
    \draw[<-] (-\x-2*\s,-2*\w)-- node[below] {$\bx^*,\blambda^*,\bz^*$} (-\x+2*\s,-2*\w);
    \draw[<-] (\x-2*\s,-2*\w)-- node[below] {$\bd^*$} (\x+2*\s,-2*\w);
    \draw[->] (-\x-2*\s,+2*\w)-- node[above] {$\bff(\cdot),\bg(\cdot)$} (-\x+2*\s,+2*\w);
    \draw[->] (\x-2*\s,+2*\w)-- node[above] {$\bM,\bp$} (\x+2*\s,+2*\w);
  \end{tikzpicture}
  \caption{Schematics of graph-based modeling and solution (top) and conventional modeling and solution (bottom).}\label{fig:schematic}
\end{figure*}
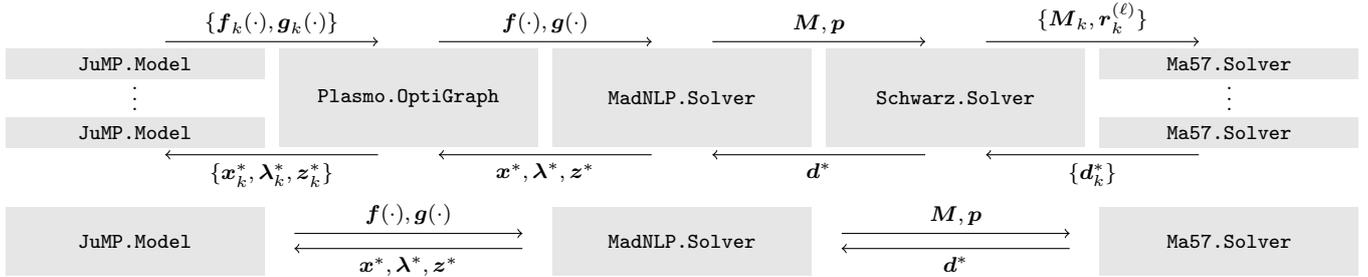

The abstraction layers within {\tt MadNLP.jl} are shown in Figure \ref{fig:schematic}. The problem is modeled as an {\tt OptiGraph} object (the core modeling object of {\tt Plasmo.jl}). The {\tt OptiGraph} is interfaced with multiple {\tt JuMP.jl} models. The {\tt JuMP.jl} model objects provide a set of local function oracles (objective, objective gradient, constraint, constraint Jacobian, and Lagrangian hessian). The {\tt Plasmo.jl}-{\tt MadNLP.jl} interface collects these local function oracles and creates a set of oracles for the full problem, where the local oracles are evaluated in parallel. 

The {\tt Solver} object of ${\tt MadNLP.jl}$ is created from the {\tt OptiGraph} object of {\tt Plasmo.jl}. The {\tt Solver} object of {\tt MadNLP.jl} uses a line-search filter IP method \citep{wachter2006implementation} to solve the problem. The step computation is performed by the linear solver specified by the user. The linear solver can be specified either as a direct solver or as the RAS solver. When the RAS solver is chosen, multiple subproblem solver objects are created by using standard direct solvers (e.g., by using {\tt Ma57} of HSL routines). These subproblem solvers are used for factorization and backsolve for $\bM_k$ blocks. The RAS scheme \eqref{eqn:ras} exploits multi-thread parallelism available in {\tt Julia}. After termination of the IP solution procedure, the primal-dual solutions are sent back to the {\tt OptiGraph} object and {\tt Model} objects from {\tt JuMP.jl} so that the user can query the solution via the interface provided by {\tt Plasmo.jl} and {\tt JuMP.jl}. See Figure \ref{fig:schematic} for a comparison with a conventional implementation.

\section{Case Studies}\label{sec:results}
The proposed computational framework was tested using transient gas network and multi-period AC power flow problems. In this section, we present the problem statements followed by numerical results and discussion.
\subsection{Transient Gas Network}\label{sec:gas}
We consider a transient gas network problem \citep{sundar2018state} of the form:
\begin{subequations}\label{eqn:gas}
\begin{align}
\min_{\substack{\rho,\varphi^a,\varphi^-,\\\alpha,s,d\in\mathbb{R}}}\; \sum_{t\in\mathcal{T}}\left(\sum_{(i,j)\in\mathcal{C}} \gamma P^a_{ijt} + \sum_{i\in\mathcal{R}}c_{it}s_{it} - \sum_{i\in D} c_{it} d_{it}\right)
\end{align}
\begin{align}
  \st\;
  &\sum_{j\in \mathcal{N}(i)}f^a_{ijt} = \sum_{j\in\mathcal{R}(i)}s_{jt} - \sum_{j\in\mathcal{D}(i)}d_{jt},i\in\mathcal{N}, t\in\mathcal{T}\\
  & \rho^{\min}_i\leq \rho_{it}\leq \rho^{\max}_i ,\;i\in\mathcal{N},\;t\in\mathcal{T}\\
  &\rho_{it}^2-\rho_{jt}^2 = -\frac{\lambda L}{D} \varphi^{a}_{ijt}|\varphi^{a}_{ijt}|,\;(i,j)\in\mathcal{P},t\in\mathcal{T}\\
  &\hat{L}(\dot{\rho}_{jt} +\dot{\rho}_{it}) = -4\varphi^-_{ijt},\;(i,j)\in\mathcal{P},\;t\in\mathcal{T}\\
  &f^a_{ijt}(\rho_{it}-\rho_{jt}) \leq 0,\; (i,j)\in \mathcal{C},\; t\in\mathcal{T}\\
  &\rho_{jt} = \alpha_{ijt}\rho_{it},\;(i,j)\in \mathcal{C},\; t\in\mathcal{T}\\
  & P^a_{ijt} \leq P_{ij}^{\max},\;(i,j)\in\mathcal{C},\;t\in\mathcal{T}\\
  &-f^a_{ij}\leq f^a_{ijt} \leq f^a_{ij},\;(i,j)\in\mathcal{C},\;t\in\mathcal{T},
\end{align}
\end{subequations}
where $\dot{\rho}_{it} = \frac{{\rho}_{it}- {\rho}_{it-1}}{\Delta t}$, $P^a_{ijt} = W_a A_{ij}$, and $f^a_{ijt} = A_{ij}\varphi^a_{ijt}$. Here, $\mathcal{N}$ is the set of junctions; $\mathcal{P}$ is the set of pipelines; $\mathcal{C}$ is the set of compressors; $\mathcal{R}$ is the set of receipts; $\mathcal{D}$ is the set of demands; $\mathcal{R}(i)$ is the set of receipts at junction $i\in\mathcal{N}$; $\mathcal{D}(i)$ is the set of demands at junction $i\in\mathcal{N}$; $\mathcal{T}$ is the time index set; $\rho$ is the densitiy; $\varphi^a$ is the average mass flux; $\varphi^-$ is the negative mass flux; $\alpha$ is the compression ratio; $s$ is the supply; $d$ is the demand; $\dot{\rho}$ is the time derivative of density; $P^a$ is the power consumption of compressor; $f$ is the mass flow; $c$ is the gas price; $\gamma$ is the economic factor; $\lambda,\hat{L},L,D,A,\Delta t$, and $W_a$ are physical parameters. To implicitly enforce the periodicity, we let $\rho_{i0}=\rho_{iT}$, where $T$ is the end time index. The gas network under study consists of 2 compressors, 6 junctions (35 junctions after discretization), 4 pipelines (32 pipelines after discretization), 1 receiving points and 5 transfer points (which work either as receipt or delivery). The model is constructed using {\tt GasModels.jl} \citep{russell2020gasmodels}. 

\subsection{Multi-Period AC  Power Flow}\label{sec:power}
We consider a multi-period AC power flow problem with storage \citep{geth2020flexible} of the form:
\begin{subequations}\label{eqn:power}
  \begin{align}
    \min_{\substack{v,s,s^g,s^s\in\mathbb{C}\\e,sc,sd,sqc\in\mathbb{R}}}\;&
      \sum_{t\in\mathcal{T}}\sum_{k\in \mathcal{G}} c^0_{kt}+c^1_{kt} \Re(s^g_{kt}) + c^2_{kt} \Re(s^g_{kt})^2\label{eqn:power-obj}
  \end{align}
  \begin{align}
    \st\;&v^L_i\leq |v_{it}|\leq v^U_i,\;i\in \mathcal{N},t\in\mathcal{T}\\
         &\hspace{-.1in}\sum_{k\in \mathcal{G}_i}s^g_{kt} -\sum_{k\in L_i}s^d_{kt} +\sum_{k\in S_i} s^s_{kt}= \sum_{j\in N_\mathcal{G}[i]} s_{ijt}, i\in \mathcal{N},t\in\mathcal{T}\\
         &\begin{aligned}
           &s_{ijt} = (Y_{ij}+Y_{ij}^c)^* \frac{|v_{it}|^2}{|T_{ij}|^2}-Y^*_{ij}\frac{v_{it}v_{jt}^*}{T_{ij}},\\
           &\qquad(i,j)\in \mathcal{E},\;t\in\mathcal{T}\\
           &s_{ijt} = (Y_{ij}+Y_{ji}^c)^* |v_{jt}|^2-Y^*_{ij}\frac{v_{it}^*v_{jt}}{T^*_{ij}},\\
           &\qquad(i,j)\in \mathcal{E}^R,\;t\in\mathcal{T}\\
         \end{aligned}\\
         &|s_{ijt}| \leq s^U_{ij},\; (i,j)\in \mathcal{E}\cup\mathcal{E}^R,\;t\in\mathcal{T}\\
         &\theta^{\Delta L}_{ij} \leq \angle (v_{it}v_{jt}^*)\leq \theta^{\Delta U}_{ij},\; (i,j)\in \mathcal{E},\;t\in\mathcal{T}\\
         &s^{gL}_k\leq s^g_{kt}\leq s^{gU}_k,\;k\in \mathcal{G},\;t\in\mathcal{T}\\
         &e_{kt}-e_{kt-1}=(\eta^c sc_t-sd_t/\eta^d)\Delta t, k\in\mathcal{S},t\in\mathcal{T}\setminus\{T\}\\
         &s^s_{kt} + (sc_{kt}-sd_{kt}) = \sqrt{-1} sqc_{kt} + s^{\text{loss}}_k,\; k\in\mathcal{S},t\in\mathcal{T}\label{eqn:power-storage}\\
         &|s^s_{kt}|\leq s^u_k,0\leq e_{kt}\leq e_k^u\;k\in\mathcal{S},t\in\mathcal{T}\\
         &0\leq sc_{kt}\leq sc_k^u,0\leq sd_{kt}\leq sc_k^u,\;k\in\mathcal{S},t\in\mathcal{T},
  \end{align}
\end{subequations}
Here, $\mathcal{G}$ is the set of generators; $\mathcal{N}$ is the set of buses; $\mathcal{E}$ is the set of (directed) branches; $\mathcal{E}^R$ is the set of branches with inverted directions; $\mathcal{S}$ is the set of storage; $\mathcal{T}$ is the time index set; $v\in\mathbb{C}$ is the voltage; $e\in\mathbb{R}$ is the state of charge; $s\in\mathbb{C}$ is the power flow; $s^g\in\mathbb{C}$ is the power generation; $s^s\in\mathbb{C}$ is the complex power injected by the storage; $sc\in\mathbb{R}$ is the charging rate; $sd\in\mathbb{R}$ is the discharging rate; $sqc\in\mathbb{R}$ is the reactive power slack; $c^0,c^1,c^2\in\mathbb{R}$ are the generation costs; $s^d\in\mathbb{C}$ is the power demand; $Y$ is the admittance; $T$ is the branch complex transformation parameter; $\eta$ is the charging efficiency; $s^{loss}$ is the storage energy loss; $\Delta t$ is the time interval. Note that \eqref{eqn:power} can be reformulated as an NLP with real variables by separately treating the real and imaginary part of the variables and equations (a polar formulation is used here). The power network under study is a variant of IEEE 14 bus test system; this comprises $14$ buses, $5$ generators, $1$ storage, $1$ shunt, and $20$ branches. The detailed model is constructed with {\tt PowerModels.jl} \citep{8442948}. 


\subsection{Results and Discussion}
We compare the proposed method ({\tt MadNLP.jl} interfaced with {\tt Plasmo.jl} and {\tt Schwarz}), with the conventional method ({\tt MadNLP.jl} interfaced with serial/parallel direct solvers {\tt Ma57} or {\tt MKL-Pardiso} along with non-graph based algebraic modeling language {\tt JuMP.jl}). The conventional methods are referred to as JuMP-Ma57 and JuMP-PardisoMKL, and the proposed method is referred to as Plasmo-Schwarz/Ma57. Furthermore, the mix of proposed/conventional approaches (JuMP-Schwarz/Ma57, Plasmo-Ma57, and Plasmo-PardisoMKL) is also tested together. For JuMP-Schwarz/Ma57, the graph partitioning tool {\tt METIS} was used to partition the primal-dual coupling graph $\widetilde{G}$ directly. A Richardson scheme was used as an iterator for the RAS scheme. The study was performed by solving the gas \eqref{eqn:gas} and power \eqref{eqn:power} problems while varying the size of the problems (by increasing the length of the prediction horizon). The code was run on a server computer equipped with 2 CPUs of Intel Xeon CPU E5-2698 v4 running 2.20GHz (20 core for each), and 20 threads are used for the computation. Code to reproduce the results can be found in \url{https://github.com/zavalab/JuliaBox/tree/master/AdchemCaseStudy}

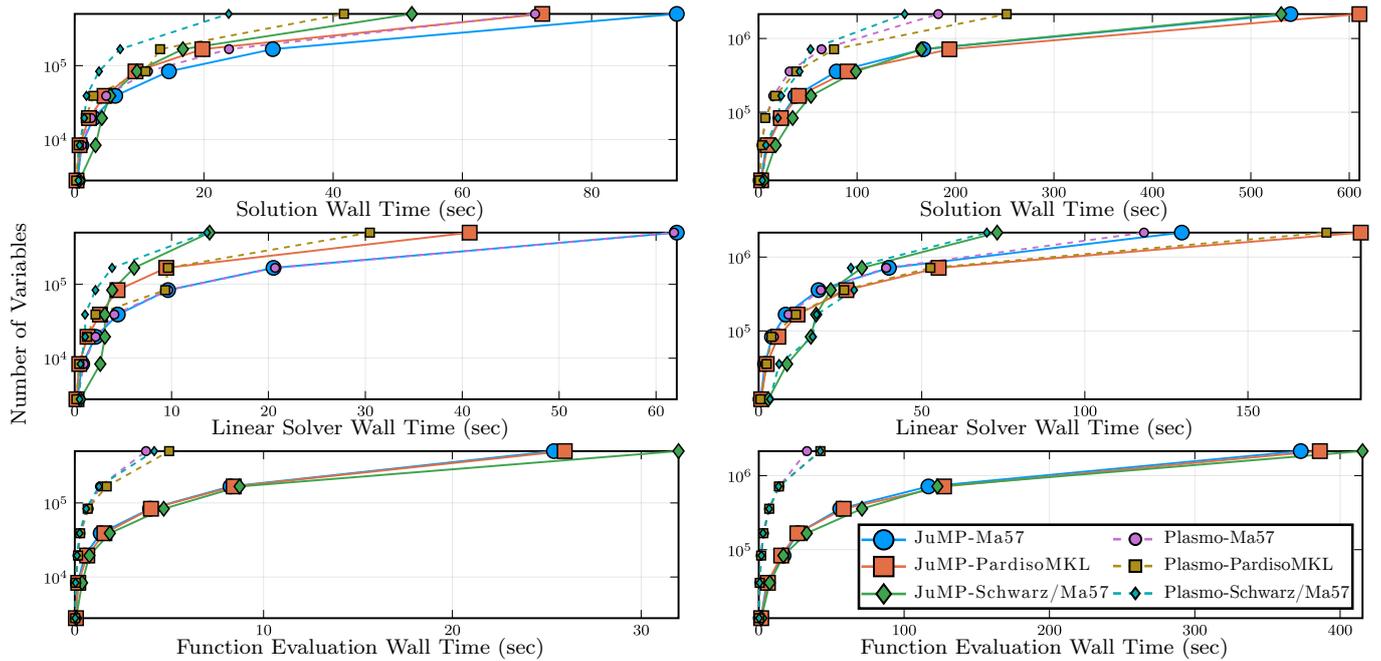
\begin{figure*}
  \begin{tikzpicture}
  \footnotesize
    \def\x{4.5};
    \def\y{2.9};
    \def\s{1.4};
    \def\lx{.8};
    \def\ly{.3};
    \node at (-\x,2*\y) {\includegraphics[width=.48\textwidth]{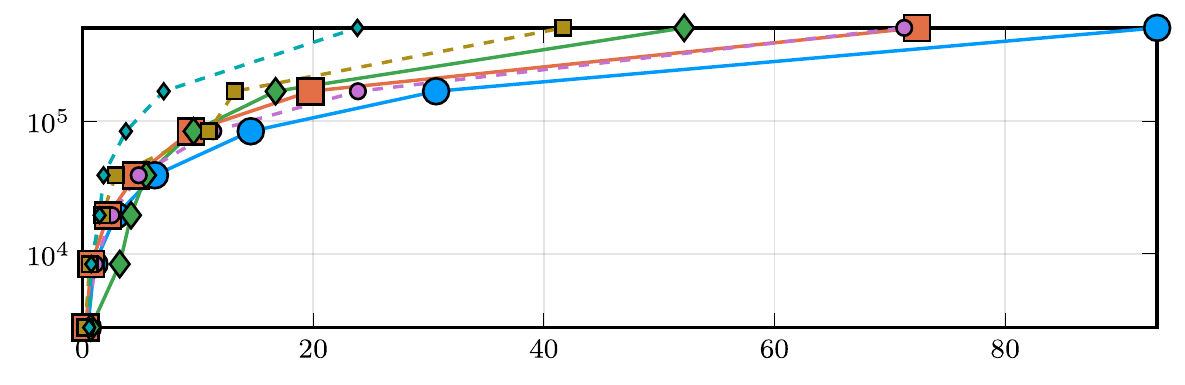}};
    \node at (\x,2*\y) {\includegraphics[width=.48\textwidth]{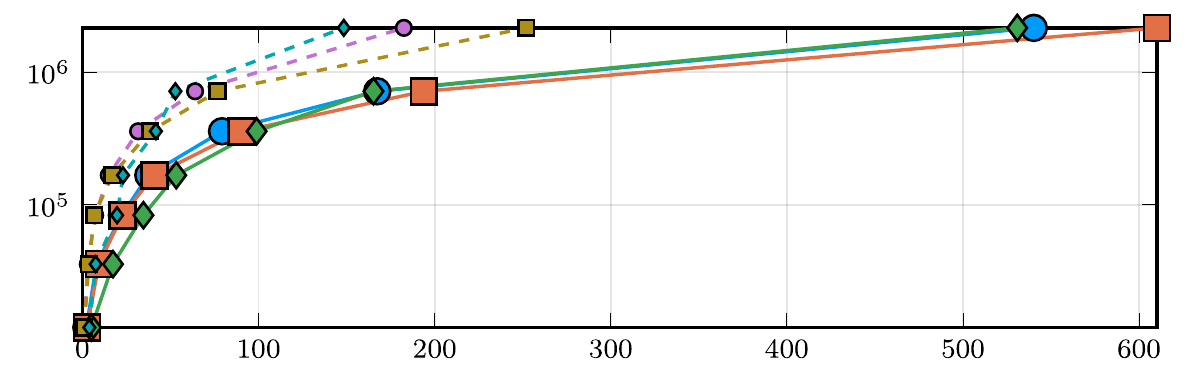}};
    \node at (-\x,\y) {\includegraphics[width=.48\textwidth]{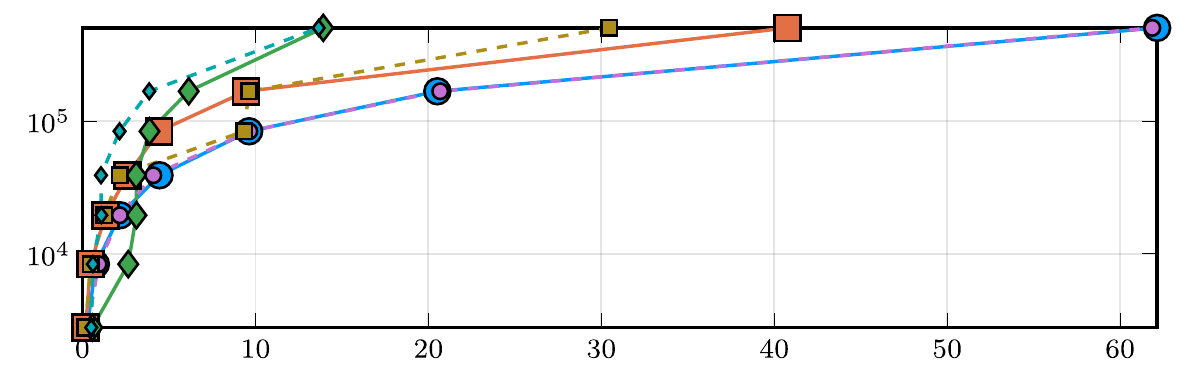}};
    \node at (\x,\y) {\includegraphics[width=.48\textwidth]{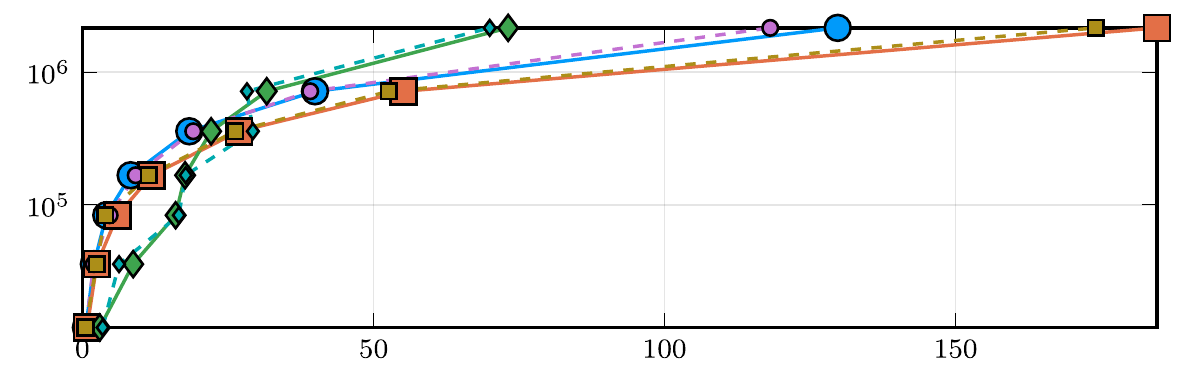}};
    \node at (-\x,0) {\includegraphics[width=.48\textwidth]{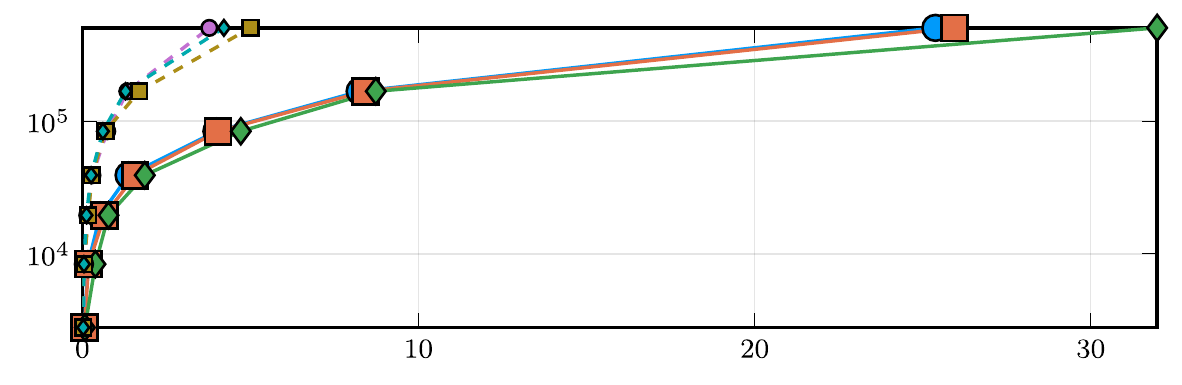}};
    \node at (\x,0) {\includegraphics[width=.48\textwidth]{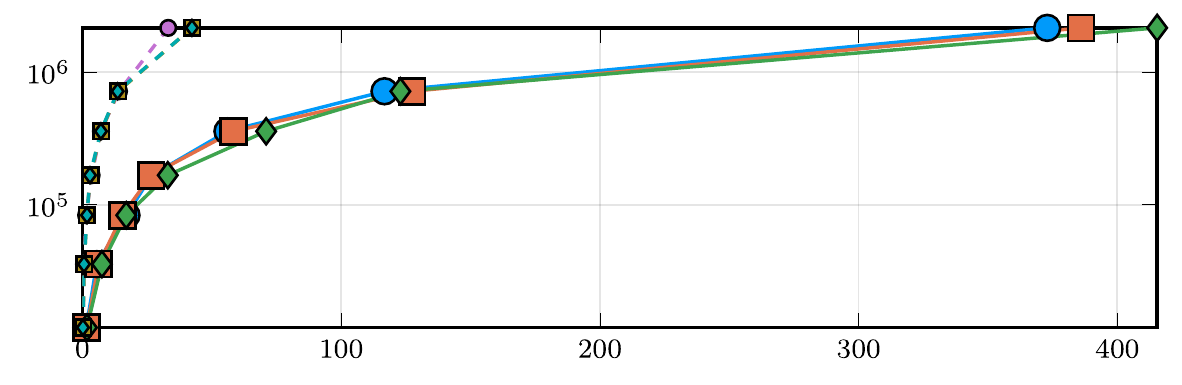}};
    \node at (-\x,0-\s) {Function Evaluation Wall Time (sec)};
    \node at (\x,0-\s) {Function Evaluation Wall Time (sec)};
    \node at (-\x,\y-\s) {Linear Solver Wall Time (sec)};
    \node at (\x,\y-\s) {Linear Solver Wall Time (sec)};
    \node at (-\x,2*\y-\s) {Solution Wall Time (sec)};
    \node at (\x,2*\y-\s) {Solution Wall Time (sec)};
    \node[rotate=90] at (-2*\x,\y) {Number of Variables};
    \node at (\x+\lx,0-\ly) {\includegraphics[width=.36\textwidth]{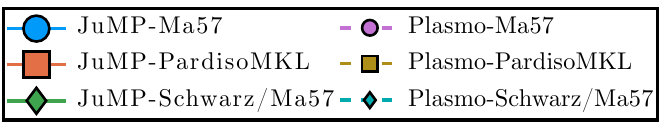}};
\end{tikzpicture}
\caption{Solution time (top), linear solver time (middle), function evaluation time (bottom) for transient gas network (left) and multi-period AC optimal power flow (right) problems.}\label{fig:result}
\end{figure*}

For both problems, we found that the graph-based approach can significantly accelerate the solution (see Figure \ref{fig:result}). In particular, comparing JuMP-Ma57 and Plasmo-Schwarz/Ma57, Plasmo-Schwarz/Ma57 becomes faster than JuMP-Ma57 when the prediction horizon is 3 days or more. Function evaluations are always faster in {\tt Plasmo.jl} compared to {\tt JuMP.jl} because the computational savings from function evaluations directly reduce the total solution time (parallelizing the function evaluation itself has no impact on the other part of the algorithm). On the other hand, one can see that the speed-up from parallel linear algebra is only observed when the problem size is sufficiently large (3 days in the gas network and 60 days in the power network). This is because the reduction in the problem size also reduces the overlap size. In our implementation, we set the size of overlap using the relative size of the block (the size of the overlap is reduced if the overall problem size is reduced). As a result, the RAS scheme \eqref{eqn:ras} may become slow, and the number of required factorization/backsolve steps increases. This indicates that the use of RAS is beneficial only when the problem size is sufficiently large. For the gas problems, the acceleration of linear algebra computations was more pronounced. In contrast, for the power problems, the acceleration of function evaluations was more pronounced. This is because the AC power flow formulation has a large number of nonlinear expressions.  By comparing the linear solver time for JuMP-Schwarz/Ma57 and Plasmo-Schwarz/Ma57, we see the advantage of using a graph-based modeling language for obtaining the partitions $\{W_k\}_{k=1}^K$. We recall that for JuMP-Schwarz/Ma57, the {\tt Metis} graph partitioning routine is directly applied to $\widetilde{G}$ while Plasmo-Schwarz/Ma57 uses the user-provided problem graph $G$. One can observe that, in general, the linear solver time is shorter for PlasmoNLP-Schwarz/Ma57. This indicates that the user-provided graph information can be leveraged for obtaining high-quality partitions.

\section{Conclusions and Future Work}\label{sec:conclusions}
We have presented a graph-based modeling and decomposition framework for large-scale nonlinear programs arising in energy infrastructures. Here, we introduce a new decomposition paradigm for linear algebra systems within an interior-point method: a restricted additive Schwarz (RAS) scheme. We implement this framework in the Julia package {\tt MadNLP.jl} and show that the RAS approach accelerates computations (compared to off-the-shelf approaches) by up to 300\%.  This work focused on applying RAS to conduct temporal decomposition; applying RAS as a spatial decomposition scheme is a future direction of interest. A surprising finding was that the RAS scheme is effective at handling instances with a large number of active inequality constraints.  We are interested in determining the reasons for this by investigating the convergence properties of the RAS scheme within an interior-point context. 

\begin{ack}
We acknowledge Andreas W\"{a}chter for helpful comments and Jordan Jalving for implementing the requested features in {\tt Plasmo.jl}. 
  
\end{ack}
\bibliography{adchem_decomp}
\end{document}